\documentclass{amsart}
\usepackage{amsfonts}
\usepackage{graphicx}
\usepackage{amscd,hyperref}

\newtheorem{theorem}{Theorem}
\theoremstyle{plain}

\newtheorem{corollary}{Corollary}

\newtheorem{lemma}{Lemma}

\newtheorem{remark}{Remark}

\numberwithin{equation}{section}

\input{tcilatex}

\begin{document}
\title[Ostrowski's Inequality]{Some Companions of Ostrowski's Inequality for
Absolutely Continuous Functions and Applications}
\author{S.S. Dragomir}
\address{School of Communications and Informatics\\
Victoria University of Technology\\
PO Box 14428, MCMC 8001\\
Victoria, Australia.}
\email{sever@matilda.vu.edu.au}
\urladdr{http://rgmia.vu.edu.au/SSDragomirWeb.html}
\date{May 22, 2002.}

\begin{abstract}
Companions of Ostrowski's integral ineqaulity for absolutely continuous
functions and applications for composite quadrature rules and for p.d.f.'s
are provided.
\end{abstract}

\maketitle

\section{Introduction}

In \cite{GS}, Guessab and Schmeisser have proved among others, the following
companion of Ostrowski's inequality.

\begin{theorem}
\label{t1.1}Let $f:\left[ a,b\right] \rightarrow \mathbb{R}$ be such that
\begin{equation}
\left| f\left( t\right) -f\left( s\right) \right| \leq M\left| t-s\right|
^{k},\;\;\;\text{for any }t,s\in \left[ a,b\right]  \label{1.1}
\end{equation}
with $k\in (0,1],$ i.e., $f\in Lip_{M}\left( k\right) .$ Then, for each $%
x\in \left[ a,\frac{a+b}{2}\right] ,$ we have the inequality
\begin{multline}
\left| \frac{f\left( x\right) +f\left( a+b-x\right) }{2}-\frac{1}{b-a}%
\int_{a}^{b}f\left( t\right) dt\right|  \label{1.2} \\
\leq \left[ \frac{2^{k+1}\left( x-a\right) ^{k+1}+\left( a+b-2x\right) ^{k+1}%
}{2^{k}\left( k+1\right) \left( b-a\right) }\right] M.
\end{multline}
This inequality is sharp for each admissable $x.$ Equality is obtained if
and only if $f=\pm Mf_{\ast }+c$ with $c\in \mathbb{R}$ and
\begin{equation}
f_{\ast }\left( t\right) =\left\{
\begin{array}{lll}
\left( x-t\right) ^{k} & \text{for} & a\leq t\leq x \\
&  &  \\
\left( t-x\right) ^{k} & \text{for} & x\leq t\leq \frac{1}{2}\left(
a+b\right) \\
&  &  \\
f_{\ast }\left( a+b-t\right) & \text{for} & \frac{1}{2}\left( a+b\right)
\leq t\leq b.%
\end{array}
\right.  \label{1.3}
\end{equation}
\end{theorem}

We remark that for $k=1,$ i.e., $f\in Lip_{M},$ since
\begin{equation*}
\frac{4\left( x-a\right) ^{2}+\left( a+b-2x\right) ^{2}}{4\left( b-a\right) }%
=\left[ \frac{1}{8}+2\left( \frac{x-\frac{3a+b}{4}}{b-a}\right) ^{2}\right]
\left( b-a\right)
\end{equation*}
then we have the inequality
\begin{multline}
\left| \frac{f\left( x\right) +f\left( a+b-x\right) }{2}-\frac{1}{b-a}%
\int_{a}^{b}f\left( t\right) dt\right|  \label{1.4} \\
\leq \left[ \frac{1}{8}+2\left( \frac{x-\frac{3a+b}{4}}{b-a}\right) ^{2}%
\right] \left( b-a\right) M
\end{multline}
for any $x\in \left[ a,\frac{a+b}{2}\right] .$

The inequality $\frac{1}{8}$ is best possible in (\ref{1.4}) in the sense
that it cannot be replaced by a smaller constant.

We must also observe that the best inequality in (\ref{1.4}) is obtained for
$x=\frac{a+3b}{4},$ giving the trapezoid type inequality
\begin{equation}
\left| \frac{f\left( \frac{3a+b}{4}\right) +f\left( \frac{a+3b}{4}\right) }{2%
}-\frac{1}{b-a}\int_{a}^{b}f\left( t\right) dt\right| \leq \frac{1}{8}\left(
b-a\right) M.  \label{1.5}
\end{equation}
The constant $\frac{1}{8}$ is sharp in (\ref{1.5}) in the sense mentioned
above.

For a recent monograph devoted to Ostrowski type inequalities, see \cite{2b}.

In this paper we improve the above results and also provide other bounds for
absolutely continuous functions whose derivatives belong to the Lebesgue
spaces $L_{p}\left[ a,b\right] ,$ $1\leq p\leq \infty .$ Some natural
applications are also provided.

\section{Some Integral Inequalities}

The following identity holds.

\begin{lemma}
\label{l2.1}Assume that $f:\left[ a,b\right] \rightarrow \mathbb{R}$ is an
absolutely continuous function on $\left[ a,b\right] $. Then we have the
equality
\begin{multline}
\frac{1}{2}\left[ f\left( x\right) +f\left( a+b-x\right) \right] -\frac{1}{%
b-a}\int_{a}^{b}f\left( t\right) dt  \label{2.1} \\
=\frac{1}{b-a}\int_{a}^{x}\left( t-a\right) f^{\prime }\left( t\right) dt+%
\frac{1}{b-a}\int_{x}^{a+b-x}\left( t-\frac{a+b}{2}\right) f^{\prime }\left(
t\right) dt \\
+\frac{1}{b-a}\int_{a+b-x}^{b}\left( t-b\right) f^{\prime }\left( t\right)
dt,
\end{multline}
for any $x\in \left[ a,\frac{a+b}{2}\right] .$
\end{lemma}

\begin{proof}
Using the integration by parts formula for Lebesgue integrals, we have
\begin{equation*}
\int_{a}^{x}\left( t-a\right) f^{\prime }\left( t\right) dt=f\left( x\right)
\left( x-a\right) -\int_{a}^{x}f\left( t\right) dt,
\end{equation*}
\begin{multline*}
\int_{x}^{a+b-x}\left( t-\frac{a+b}{2}\right) f^{\prime }\left( t\right) dt
\\
=f\left( a+b-x\right) \left( \frac{a+b}{2}-x\right) -f\left( x\right) \left(
x-\frac{a+b}{2}\right) -\int_{x}^{a+b-x}f\left( t\right) dt
\end{multline*}
and
\begin{equation*}
\int_{a+b-x}^{b}\left( t-b\right) f^{\prime }\left( t\right) dt=\left(
x-a\right) f\left( a+b-x\right) -\int_{a+b-x}^{b}f\left( t\right) dt.
\end{equation*}
Summing the above equalities, we deduce the desired identity (\ref{2.1}).
\end{proof}

\begin{remark}
\label{r2.2}The identity (\ref{2.1}) was obtained in \cite[Lemma 3.2]{GS}
for the case of piecewise continuously differentiable functions on $\left[
a,b\right] .$
\end{remark}

The following result holds.

\begin{theorem}
\label{t2.3}Let $f:\left[ a,b\right] \rightarrow \mathbb{R}$ be an
absolutely continuous function on $\left[ a,b\right] .$ Then we have the
inequality
\begin{align}
& \left| \frac{1}{2}\left[ f\left( x\right) +f\left( a+b-x\right) \right] -%
\frac{1}{b-a}\int_{a}^{b}f\left( t\right) dt\right|  \label{2.2} \\
& \leq \frac{1}{b-a}\left[ \int_{a}^{x}\left( t-a\right) \left| f^{\prime
}\left( t\right) \right| dt\right.  \notag \\
& \;\;\;\;\;\;\;\;\;\;\;+\left. \int_{x}^{a+b-x}\left| t-\frac{a+b}{2}%
\right| \left| f^{\prime }\left( t\right) \right| dt+\int_{a+b-x}^{b}\left(
b-t\right) \left| f^{\prime }\left( t\right) \right| dt\right]  \notag \\
& :=M\left( x\right)  \notag
\end{align}
for any $x\in \left[ a,\frac{a+b}{2}\right] .$

If $f^{\prime }\in L_{\infty }\left[ a,b\right] ,$ then we have the
inequalities
\begin{multline}
M\left( x\right) \leq \frac{1}{b-a}\left[ \frac{\left( x-a\right) ^{2}}{2}%
\left\| f^{\prime }\right\| _{\left[ a,x\right] ,\infty }+\left( \frac{a+b}{2%
}-x\right) ^{2}\left\| f^{\prime }\right\| _{\left[ x,a+b-x\right] ,\infty
}\right.  \label{2.3} \\
+\left. \frac{\left( x-a\right) ^{2}}{2}\left\| f^{\prime }\right\| _{\left[
a+b-x,b\right] ,\infty }\right]
\end{multline}
\begin{equation*}
\leq \left\{
\begin{array}{l}
\left[ \dfrac{1}{8}+2\left( \dfrac{x-\frac{3a+b}{4}}{b-a}\right) ^{2}\right]
\left( b-a\right) \left\| f^{\prime }\right\| _{\left[ a,b\right] ,\infty }
\\
\\
\left[ \dfrac{1}{2^{\alpha -1}}\left( \dfrac{x-a}{b-a}\right) ^{2\alpha
}+\left( \dfrac{x-\frac{a+b}{2}}{b-a}\right) ^{2\alpha }\right] ^{\frac{1}{%
\alpha }} \\
\times \left[ \left\| f^{\prime }\right\| _{\left[ a,x\right] ,\infty
}^{\beta }+\left\| f^{\prime }\right\| _{\left[ x,a+b-x\right] ,\infty
}^{\beta }+\left\| f^{\prime }\right\| _{\left[ a+b-x,b\right] ,\infty
}^{\beta }\right] ^{\frac{1}{\beta }}\left( b-a\right) \\
\hfill \text{if \hspace{0.05in}}\alpha >1,\;\frac{1}{\alpha }+\frac{1}{\beta
}=1, \\
\\
\max \left\{ \dfrac{1}{2}\left( \dfrac{x-a}{b-a}\right) ^{2},\left( \dfrac{x-%
\frac{a+b}{2}}{b-a}\right) ^{2}\right\} \\
\hfill \times \left[ \left\| f^{\prime }\right\| _{\left[ a,x\right] ,\infty
}+\left\| f^{\prime }\right\| _{\left[ x,a+b-x\right] ,\infty }+\left\|
f^{\prime }\right\| _{\left[ a+b-x,b\right] ,\infty }\right] \left(
b-a\right)%
\end{array}
\right.
\end{equation*}
for any $x\in \left[ a,\frac{a+b}{2}\right] .$

The inequality (\ref{2.2}), the first inequality in (\ref{2.3}) and the
constant $\frac{1}{8}$ are sharp.
\end{theorem}

\begin{proof}
The inequality (\ref{2.2}) follows by Lemma \ref{l2.1} on taking the modulus
and using it properties.

If $f^{\prime }\in L_{\infty }\left[ a,b\right] ,$ then
\begin{equation*}
\int_{a}^{x}\left( t-a\right) \left| f^{\prime }\left( t\right) \right|
dt\leq \frac{\left( x-a\right) ^{2}}{2}\left\| f^{\prime }\right\| _{\left[
a,x\right] ,\infty },
\end{equation*}
\begin{equation*}
\int_{x}^{a+b-x}\left| t-\frac{a+b}{2}\right| \left| f^{\prime }\left(
t\right) \right| dt\leq \left( \frac{a+b}{2}-x\right) ^{2}\left\| f^{\prime
}\right\| _{\left[ x,a+b-x\right] ,\infty },
\end{equation*}
\begin{equation*}
\int_{a+b-x}^{b}\left( b-t\right) \left| f^{\prime }\left( t\right) \right|
dt\leq \frac{\left( x-a\right) ^{2}}{2}\left\| f^{\prime }\right\| _{\left[
a+b-x,b\right] ,\infty }
\end{equation*}
and the first inequality in (\ref{2.3}) is proved.

Denote
\begin{multline*}
\tilde{M}\left( x\right) :=\frac{\left( x-a\right) ^{2}}{2}\left\| f^{\prime
}\right\| _{\left[ a,x\right] ,\infty }+\left( \frac{a+b}{2}-x\right)
^{2}\left\| f^{\prime }\right\| _{\left[ x,a+b-x\right] ,\infty } \\
+\frac{\left( x-a\right) ^{2}}{2}\left\| f^{\prime }\right\| _{\left[ a+b-x,b%
\right] ,\infty }
\end{multline*}
for $x\in \left[ a,\frac{a+b}{2}\right] .$

Firstly, observe that
\begin{eqnarray*}
\tilde{M}\left( x\right) &\leq &\max \left\{ \left\| f^{\prime }\right\| _{%
\left[ a,x\right] ,\infty },\left\| f^{\prime }\right\| _{\left[ x,a+b-x%
\right] ,\infty },\left\| f^{\prime }\right\| _{\left[ a+b-x,b\right]
,\infty }\right\} \\
&&\times \left[ \frac{\left( x-a\right) ^{2}}{2}+\left( \frac{a+b}{2}%
-x\right) ^{2}+\frac{\left( x-a\right) ^{2}}{2}\right] \\
&=&\left\| f^{\prime }\right\| _{\left[ a,b\right] ,\infty }\left[ \dfrac{1}{%
8}\left( b-a\right) ^{2}+2\left( x-\dfrac{3a+b}{4}\right) ^{2}\right]
\end{eqnarray*}
and the first inequality in (\ref{2.3}) is proved.

Using H\"{o}lder's inequality for $\alpha >1,$ $\frac{1}{\alpha }+\frac{1}{%
\beta }=1,$ we also have
\begin{multline*}
\tilde{M}\left( x\right) \leq \left\{ \left[ \frac{\left( x-a\right) ^{2}}{2}%
\right] ^{\alpha }+\left( x-\frac{a+b}{2}\right) ^{2\alpha }+\left[ \frac{%
\left( x-a\right) ^{2}}{2}\right] ^{\alpha }\right\} ^{\frac{1}{\alpha }} \\
\times \left[ \left\| f^{\prime }\right\| _{\left[ a,x\right] ,\infty
}^{\beta }+\left\| f^{\prime }\right\| _{\left[ x,a+b-x\right] ,\infty
}^{\beta }+\left\| f^{\prime }\right\| _{\left[ a+b-x,b\right] ,\infty
}^{\beta }\right] ^{\frac{1}{\beta }}
\end{multline*}
giving the second inequality in (\ref{2.3}).

Finally, we also observe that
\begin{multline*}
\tilde{M}\left( x\right) \leq \max \left\{ \frac{\left( x-a\right) ^{2}}{2}%
,\left( x-\frac{a+b}{2}\right) ^{2}\right\} \\
\times \left[ \left\| f^{\prime }\right\| _{\left[ a,x\right] ,\infty
}+\left\| f^{\prime }\right\| _{\left[ x,a+b-x\right] ,\infty }+\left\|
f^{\prime }\right\| _{\left[ a+b-x,b\right] ,\infty }\right] .
\end{multline*}
The sharpness of the inequalities mentioned follows from Theorem \ref{t1.1}
for $k=1.$ We omit the details.
\end{proof}

\begin{remark}
\label{r2.4}If in Theorem \ref{t2.3} we choose $x=a,$ then we get
\begin{equation}
\left\vert \frac{f\left( a\right) +f\left( b\right) }{2}-\frac{1}{b-a}%
\int_{a}^{b}f\left( t\right) dt\right\vert \leq \frac{1}{4}\left( b-a\right)
\left\Vert f^{\prime }\right\Vert _{\left[ a,b\right] ,\infty }  \label{2.4}
\end{equation}%
with $\frac{1}{4}$ as a sharp constant (see for example \cite[p. 25]{2b}).

If in the same theorem we now choose $x=\frac{a+b}{2},$ then we get
\begin{align}
\left\vert f\left( \frac{a+b}{2}\right) -\frac{1}{b-a}\int_{a}^{b}f\left(
t\right) dt\right\vert & \leq \frac{1}{8}\left( b-a\right) \left[ \left\Vert
f^{\prime }\right\Vert _{\left[ a,\frac{a+b}{2}\right] ,\infty }+\left\Vert
f^{\prime }\right\Vert _{\left[ \frac{a+b}{2},b\right] ,\infty }\right]
\label{2.5} \\
& \leq \frac{1}{4}\left( b-a\right) \left\Vert f^{\prime }\right\Vert _{%
\left[ a,b\right] ,\infty }  \notag
\end{align}%
with the constants $\frac{1}{8}$ and $\frac{1}{4}$ being sharp. This result
was obtained in \cite{SSD2}.
\end{remark}

It is natural to consider the following corollary.

\begin{corollary}
\label{c2.4.a}With the assumptions in Theorem \ref{t2.3}, one has the
inequality:
\begin{equation}
\left| \frac{f\left( \frac{3a+b}{4}\right) +f\left( \frac{a+3b}{4}\right) }{2%
}-\frac{1}{b-a}\int_{a}^{b}f\left( t\right) dt\right| \leq \frac{1}{8}\left(
b-a\right) \left\| f^{\prime }\right\| _{\left[ a,b\right] ,\infty }.
\label{2.5.a}
\end{equation}
The constant $\frac{1}{8}$ is best possible in the sense that it cannot be
replaced by a smaller constant.
\end{corollary}

The case when $f^{\prime }\in L_{p}\left[ a,b\right] $ , $p>1$ is embodied
in the following theorem.

\begin{theorem}
\label{t2.5}Let $f:\left[ a,b\right] \rightarrow \mathbb{R}$ be an
absolutely continuous function on $\left[ a,b\right] $ so that $f^{\prime
}\in L_{p}\left[ a,b\right] ,$ $p>1.$ If $M\left( x\right) $ is as defined
in (\ref{2.2}), then we have the bounds:
\begin{multline}
M\left( x\right) \leq \left. \frac{1}{\left( q+1\right) ^{\frac{1}{q}}}%
\right[ \left( \frac{x-a}{b-a}\right) ^{1+\frac{1}{q}}\left\| f^{\prime
}\right\| _{\left[ a,x\right] ,p}  \label{2.6} \\
\left. +2^{\frac{1}{q}}\left( \frac{\frac{a+b}{2}-x}{b-a}\right) ^{1+\frac{1%
}{q}}\left\| f^{\prime }\right\| _{\left[ x,a+b-x\right] ,p}+\left( \frac{x-a%
}{b-a}\right) ^{1+\frac{1}{q}}\left\| f^{\prime }\right\| _{\left[ a+b-x,b%
\right] ,p}\right] \left( b-a\right) ^{\frac{1}{q}}
\end{multline}
\begin{equation*}
\leq \frac{1}{\left( q+1\right) ^{\frac{1}{q}}}\left\{
\begin{array}{l}
\left[ 2\left( \frac{x-a}{b-a}\right) ^{1+\frac{1}{q}}+2^{\frac{1}{q}}\left(
\frac{\frac{a+b}{2}-x}{b-a}\right) ^{1+\frac{1}{q}}\right] \\
\hfill \times \max \left\{ \left\| f^{\prime }\right\| _{\left[ a,x\right]
,p},\left\| f^{\prime }\right\| _{\left[ x,a+b-x\right] ,p},\left\|
f^{\prime }\right\| _{\left[ a+b-x,b\right] ,p}\right\} \left( b-a\right) ^{%
\frac{1}{q}} \\
\\
\left[ 2\left( \frac{x-a}{b-a}\right) ^{\alpha +\frac{\alpha }{q}}+2^{\frac{%
\alpha }{q}}\left( \frac{\frac{a+b}{2}-x}{b-a}\right) ^{\alpha +\frac{\alpha
}{q}}\right] ^{\frac{1}{\alpha }} \\
\times \left[ \left\| f^{\prime }\right\| _{\left[ a,x\right] ,p}^{\beta
}+\left\| f^{\prime }\right\| _{\left[ x,a+b-x\right] ,p}^{\beta }+\left\|
f^{\prime }\right\| _{\left[ a+b-x,b\right] ,p}^{\beta }\right] ^{\frac{1}{%
\beta }}\left( b-a\right) ^{\frac{1}{q}} \\
\hfill \text{if \hspace{0.05in}}\alpha >1,\;\frac{1}{\alpha }+\frac{1}{\beta
}=1, \\
\\
\max \left\{ \left( \frac{x-a}{b-a}\right) ^{1+\frac{1}{q}},2^{\frac{1}{q}%
}\left( \frac{\frac{a+b}{2}-x}{b-a}\right) ^{1+\frac{1}{q}}\right\} \\
\hfill \times \left[ \left\| f^{\prime }\right\| _{\left[ a,x\right]
,p}+\left\| f^{\prime }\right\| _{\left[ x,a+b-x\right] ,p}+\left\|
f^{\prime }\right\| _{\left[ a+b-x,b\right] ,p}\right] \left( b-a\right) ^{%
\frac{1}{q}}%
\end{array}
\right.
\end{equation*}
for any $x\in \left[ a,\frac{a+b}{2}\right] .$
\end{theorem}

\begin{proof}
Using H\"{o}lder's integral inequality for $p>1,$ $\frac{1}{p}+\frac{1}{q}%
=1, $ we have
\begin{equation*}
\int_{a}^{x}\left( t-a\right) \left| f^{\prime }\left( t\right) \right|
dt\leq \left( \int_{a}^{x}\left( t-a\right) ^{q}dt\right) ^{\frac{1}{q}%
}\left\| f^{\prime }\right\| _{\left[ a,x\right] ,p}=\frac{\left( x-a\right)
^{1+\frac{1}{q}}}{\left( q+1\right) ^{\frac{1}{q}}}\left\| f^{\prime
}\right\| _{\left[ a,x\right] ,p},
\end{equation*}
\begin{align*}
\int_{x}^{a+b-x}\left| t-\frac{a+b}{2}\right| \left| f^{\prime }\left(
t\right) \right| dt& \leq \left( \int_{x}^{a+b-x}\left| t-\frac{a+b}{2}%
\right| ^{q}dt\right) ^{\frac{1}{q}}\left\| f^{\prime }\right\| _{\left[
x,a+b-x\right] ,p} \\
& =\frac{2^{\frac{1}{q}}\left( \frac{a+b}{2}-x\right) ^{1+\frac{1}{q}}}{%
\left( q+1\right) ^{\frac{1}{q}}}\left\| f^{\prime }\right\| _{\left[ x,a+b-x%
\right] ,p}
\end{align*}
and
\begin{align*}
\int_{a+b-x}^{b}\left( b-t\right) \left| f^{\prime }\left( t\right) \right|
dt& \leq \left( \int_{a+b-x}^{b}\left( b-t\right) ^{q}dt\right) ^{\frac{1}{q}%
}\left\| f^{\prime }\right\| _{\left[ a+b-x,b\right] ,p} \\
& =\frac{\left( x-a\right) ^{1+\frac{1}{q}}}{\left( q+1\right) ^{\frac{1}{q}}%
}\left\| f^{\prime }\right\| _{\left[ a+b-x,b\right] ,p}.
\end{align*}
Summing the above inequalities, we deduce the first bound in (\ref{2.6}).

The last part may be proved in a similar fashion to the one in Theorem \ref%
{t2.3}, and we omit the details.
\end{proof}

\begin{remark}
\label{r2.6}If in (\ref{2.6}) we choose $\alpha =q,$ $\beta =p,$ $\frac{1}{p}%
+\frac{1}{q}=1,$ $p>1,$ then we get the inequality
\begin{equation}
M\left( x\right) \leq \frac{2^{\frac{1}{q}}}{\left( q+1\right) ^{\frac{1}{q}}%
}\left[ \left( \dfrac{x-a}{b-a}\right) ^{q+1}+\left( \dfrac{\frac{a+b}{2}-x}{%
b-a}\right) ^{q+1}\right] ^{\frac{1}{q}}\left( b-a\right) ^{\frac{1}{q}%
}\left\| f^{\prime }\right\| _{\left[ a,b\right] ,p}  \label{2.7}
\end{equation}
for any $x\in \left[ a,\frac{a+b}{2}\right] .$
\end{remark}

\begin{remark}
\label{r2.7}If in Theorem \ref{t2.5} we choose $x=a,$ then we get the
trapezoid inequality
\begin{equation}
\left\vert \frac{f\left( a\right) +f\left( b\right) }{2}-\frac{1}{b-a}%
\int_{a}^{b}f\left( t\right) dt\right\vert \leq \frac{1}{2}\cdot \frac{%
\left( b-a\right) ^{\frac{1}{q}}\left\Vert f^{\prime }\right\Vert _{\left[
a,b\right] ,p}}{\left( q+1\right) ^{\frac{1}{q}}},  \label{2.8}
\end{equation}%
The constant $\frac{1}{2}$ is best possible in the sense that it cannot be
replaced by a smaller constant (see for example \cite[p. 42]{2b}).
\end{remark}

Indeed, if we assume that (\ref{2.8}) holds with a constant $C>0,$ instead
of $\frac{1}{2},$ i.e.,
\begin{equation}
\left| \frac{f\left( a\right) +f\left( b\right) }{2}-\frac{1}{b-a}%
\int_{a}^{b}f\left( t\right) dt\right| \leq C\cdot \frac{\left( b-a\right) ^{%
\frac{1}{q}}\left\| f^{\prime }\right\| _{\left[ a,b\right] ,p}}{\left(
q+1\right) ^{\frac{1}{q}}},  \label{2.9}
\end{equation}
then for the function $f:\left[ a,b\right] \rightarrow \mathbb{R}$, $f\left(
x\right) =k\left| x-\frac{a+b}{2}\right| ,$ $k>0,$ we have
\begin{align*}
\frac{f\left( a\right) +f\left( b\right) }{2}& =k\cdot \frac{b-a}{2}, \\
\frac{1}{b-a}\int_{a}^{b}f\left( t\right) dt& =k\cdot \frac{b-a}{4}, \\
\left\| f^{\prime }\right\| _{\left[ a,b\right] ,p}& =k\left( b-a\right) ^{%
\frac{1}{p}};
\end{align*}
and by (\ref{2.9}) we deduce
\begin{equation*}
\left| \frac{k\left( b-a\right) }{2}-\frac{k\left( b-a\right) }{4}\right|
\leq \frac{C\cdot k\left( b-a\right) }{\left( q+1\right) ^{\frac{1}{q}}},
\end{equation*}
giving $C\geq \frac{\left( q+1\right) ^{\frac{1}{q}}}{4}.$ Letting $%
q\rightarrow 1+,$ we deduce $C\geq \frac{1}{2},$ and the sharpness of the
constant is proved.

\begin{remark}
\label{r2.8}If in Theorem \ref{t2.5} we choose $x=\frac{a+b}{2},$ then we
get the midpoint inequality
\begin{align}
& \left| f\left( \frac{a+b}{2}\right) -\frac{1}{b-a}\int_{a}^{b}f\left(
t\right) dt\right|  \label{2.10} \\
& \leq \frac{1}{2}\cdot \frac{\left( b-a\right) ^{\frac{1}{q}}}{2^{\frac{1}{q%
}}\left( q+1\right) ^{\frac{1}{q}}}\left[ \left\| f^{\prime }\right\| _{%
\left[ a,\frac{a+b}{2}\right] ,p}+\left\| f^{\prime }\right\| _{\left[ \frac{%
a+b}{2},b\right] ,p}\right]  \notag \\
& \leq \frac{1}{2}\cdot \frac{\left( b-a\right) ^{\frac{1}{q}}}{\left(
q+1\right) ^{\frac{1}{q}}}\left\| f^{\prime }\right\| _{\left[ a,b\right]
,p},\;\;p>1,\;\frac{1}{p}+\frac{1}{q}=1.  \notag
\end{align}
\end{remark}

In both inequalities the constant $\frac{1}{2}$ is sharp in the sense that
it cannot be replaced by a smaller constant.

To show this fact, assume that (\ref{2.10}) holds with $C,D>0,$ i.e.,
\begin{align}
& \left| f\left( \frac{a+b}{2}\right) -\frac{1}{b-a}\int_{a}^{b}f\left(
t\right) dt\right|  \label{2.11} \\
& \leq C\cdot \frac{\left( b-a\right) ^{\frac{1}{q}}}{2^{\frac{1}{q}}\left(
q+1\right) ^{\frac{1}{q}}}\left[ \left\| f^{\prime }\right\| _{\left[ a,%
\frac{a+b}{2}\right] ,p}+\left\| f^{\prime }\right\| _{\left[ \frac{a+b}{2},b%
\right] ,p}\right]  \notag \\
& \leq D\cdot \frac{\left( b-a\right) ^{\frac{1}{q}}}{\left( q+1\right) ^{%
\frac{1}{q}}}\left\| f^{\prime }\right\| _{\left[ a,b\right] ,p}.  \notag
\end{align}
For the function $f:\left[ a,b\right] \rightarrow \mathbb{R}$, $f\left(
x\right) =k\left| x-\frac{a+b}{2}\right| ,$ $k>0,$ we have
\begin{equation*}
f\left( \frac{a+b}{2}\right) =0,\;\;\frac{1}{b-a}\int_{a}^{b}f\left(
t\right) dt=\frac{k\left( b-a\right) }{4},
\end{equation*}
\begin{equation*}
\left\| f^{\prime }\right\| _{\left[ a,\frac{a+b}{2}\right] ,p}+\left\|
f^{\prime }\right\| _{\left[ \frac{a+b}{2},b\right] ,p}=2\left( \frac{b-a}{2}%
\right) ^{\frac{1}{p}}k=2^{\frac{1}{q}}\left( b-a\right) ^{\frac{1}{p}}k,
\end{equation*}
\begin{equation*}
\left\| f^{\prime }\right\| _{\left[ a,b\right] ,p}=\left( b-a\right) ^{%
\frac{1}{p}}k;
\end{equation*}
and then by (\ref{2.11}) we deduce
\begin{equation*}
\frac{k\left( b-a\right) }{4}\leq C\cdot \frac{k\left( b-a\right) }{\left(
q+1\right) ^{\frac{1}{q}}}\leq D\cdot \frac{k\left( b-a\right) }{\left(
q+1\right) ^{\frac{1}{q}}},
\end{equation*}
giving $C,D\geq \frac{\left( q+1\right) ^{\frac{1}{q}}}{4}$ for any $q>1.$
Letting $q\rightarrow 1+,$ we deduce $C,D\geq \frac{1}{2}$ and the sharpness
of the constants in (\ref{2.10}) are proved.

The following result is useful in providing the best quadrature rule in the
class for approximating the integral of an absolutely continuous function
whose derivative is in $L_{p}\left[ a,b\right] .$

\begin{corollary}
\label{c2.8}Assume that $f:\left[ a,b\right] \rightarrow \mathbb{R}$ is an
absolutely continuous function so that $f^{\prime }\in L_{p}\left[ a,b\right]
,$ $p>1.$ Then one has the inequality
\begin{equation}
\left| \frac{f\left( \frac{3a+b}{4}\right) +f\left( \frac{a+3b}{4}\right) }{2%
}-\frac{1}{b-a}\int_{a}^{b}f\left( t\right) dt\right| \leq \frac{1}{4}\frac{%
\left( b-a\right) ^{\frac{1}{q}}}{\left( q+1\right) ^{\frac{1}{q}}}\left\|
f^{\prime }\right\| _{\left[ a,b\right] ,p},  \label{2.12}
\end{equation}
where $\frac{1}{p}+\frac{1}{q}=1.$

The constant $\frac{1}{4}$ is the best possible in the sense that it cannot
be replaced by a smaller constant.
\end{corollary}

\begin{proof}
The inequality follows by Theorem \ref{t2.5} and Remark \ref{r2.6} on
choosing $x=\frac{3a+b}{4}.$

To prove the sharpness of the constant, assume that (\ref{2.12}) holds with
a constant $E>0,$ i.e.,
\begin{equation}
\left| \frac{f\left( \frac{3a+b}{4}\right) +f\left( \frac{a+3b}{4}\right) }{2%
}-\frac{1}{b-a}\int_{a}^{b}f\left( t\right) dt\right| \leq E\cdot \frac{%
\left( b-a\right) ^{\frac{1}{q}}}{\left( q+1\right) ^{\frac{1}{q}}}\left\|
f^{\prime }\right\| _{\left[ a,b\right] ,p}.  \label{2.13}
\end{equation}
Consider the function $f:\left[ a,b\right] \rightarrow \mathbb{R}$,
\begin{equation*}
f\left( x\right) =\left\{
\begin{array}{lll}
\left| x-\dfrac{3a+b}{4}\right| & \text{if} & x\in \left[ a,\frac{a+b}{2}%
\right] \\
&  &  \\
\left| x-\dfrac{a+3b}{4}\right| & \text{if} & x\in \left( \frac{a+b}{2},b%
\right] .%
\end{array}
\right.
\end{equation*}
Then $f$ is absolutely continuous and $f^{\prime }\in L_{p}\left[ a,b\right]
,$ $p>1.$ We also have
\begin{equation*}
\frac{1}{2}\left[ f\left( \frac{3a+b}{4}\right) +f\left( \frac{a+3b}{4}%
\right) \right] =0,\;\;\frac{1}{b-a}\int_{a}^{b}f\left( t\right) dt=\frac{b-a%
}{8}
\end{equation*}
\begin{equation*}
\left\| f^{\prime }\right\| _{\left[ a,b\right] ,p}=\left( b-a\right) ^{%
\frac{1}{p}},
\end{equation*}
and then, by (\ref{2.13}), we obtain:
\begin{equation*}
\frac{b-a}{8}\leq E\frac{\left( b-a\right) }{\left( q+1\right) ^{\frac{1}{q}}%
}
\end{equation*}
giving $E\geq \frac{\left( q+1\right) ^{\frac{1}{q}}}{8}$ for any $q>1,$
i.e., $E\geq \frac{1}{4},$ and the corollary is proved.
\end{proof}

If one is interested in obtaining bounds in terms of the $1-$norm for the
derivative, then the following result may be useful.

\begin{theorem}
\label{t2.9}Assume that the function $f:\left[ a,b\right] \rightarrow
\mathbb{R}$ is absolutely continuous on $\left[ a,b\right] .$ If $M\left(
x\right) $ is as in equation (\ref{2.2}), then we have the bounds
\begin{multline}
M\left( x\right) \leq \left( \frac{x-a}{b-a}\right) \left\| f^{\prime
}\right\| _{\left[ a,x\right] ,1}  \label{2.14} \\
+\left( \frac{\frac{a+b}{2}-x}{b-a}\right) \left\| f^{\prime }\right\| _{%
\left[ x,a+b-x\right] ,1}+\left( \frac{x-a}{b-a}\right) \left\| f^{\prime
}\right\| _{\left[ a+b-x,b\right] ,1}
\end{multline}
\begin{equation*}
\leq \left\{
\begin{array}{l}
\left[ \dfrac{1}{4}+\left| \dfrac{x-\frac{3a+b}{4}}{b-a}\right| \right]
\left\| f^{\prime }\right\| _{\left[ a,b\right] ,1} \\
\\
\left[ 2\left( \dfrac{x-a}{b-a}\right) ^{\alpha }+\left( \dfrac{\frac{a+b}{2}%
-x}{b-a}\right) ^{\alpha }\right] ^{\frac{1}{\alpha }} \\
\times \left[ \left\| f^{\prime }\right\| _{\left[ a,x\right] ,1}^{\beta
}+\left\| f^{\prime }\right\| _{\left[ x,a+b-x\right] ,1}^{\beta }+\left\|
f^{\prime }\right\| _{\left[ a+b-x,b\right] ,1}^{\beta }\right] ^{\frac{1}{%
\beta }} \\
\hfill \text{if \hspace{0.05in}}\alpha >1,\;\frac{1}{\alpha }+\frac{1}{\beta
}=1, \\
\\
\left[ \dfrac{x+\frac{b-3a}{2}}{b-a}\right] \max \left[ \left\| f^{\prime
}\right\| _{\left[ a,x\right] ,1},\left\| f^{\prime }\right\| _{\left[
x,a+b-x\right] ,1},\left\| f^{\prime }\right\| _{\left[ a+b-x,b\right] ,1}%
\right] .%
\end{array}
\right.
\end{equation*}
\end{theorem}

The proof is as in Theorem \ref{t2.3} and we omit it.

\begin{remark}
\label{r2.10}By the use of Theorem \ref{t2.5}, for $x=a,$ we get the
trapezoid inequality (see for example \cite[p. 55]{2b})
\begin{equation}
\left\vert \frac{f\left( a\right) +f\left( b\right) }{2}-\frac{1}{b-a}%
\int_{a}^{b}f\left( t\right) dt\right\vert \leq \frac{1}{2}\left\Vert
f^{\prime }\right\Vert _{\left[ a,b\right] ,1}.  \label{2.15}
\end{equation}%
If in (\ref{2.14}) we also choose $x=\frac{a+b}{2},$ then we get the mid
point inequality (see for example \cite[p. 56]{2b})
\begin{equation}
\left\vert f\left( \frac{a+b}{2}\right) -\frac{1}{b-a}\int_{a}^{b}f\left(
t\right) dt\right\vert \leq \frac{1}{2}\left\Vert f^{\prime }\right\Vert _{%
\left[ a,b\right] ,1}.  \label{2.16}
\end{equation}
\end{remark}

The following corollary also holds.

\begin{corollary}
\label{c2.11}With the assumption in Theorem \ref{t2.5}, one has the
inequality:
\begin{equation}
\left| \frac{f\left( \frac{3a+b}{4}\right) +f\left( \frac{a+3b}{4}\right) }{2%
}-\frac{1}{b-a}\int_{a}^{b}f\left( t\right) dt\right| \leq \frac{1}{4}%
\left\| f^{\prime }\right\| _{\left[ a,b\right] ,1}.  \label{2.17}
\end{equation}
\end{corollary}

\section{A Composite Quadrature Formula}

We use the following inequalities obtained in the previous section:
\begin{multline}
\left| \frac{f\left( \frac{3a+b}{4}\right) +f\left( \frac{a+3b}{4}\right) }{2%
}-\frac{1}{b-a}\int_{a}^{b}f\left( t\right) dt\right|  \label{3.1} \\
\leq \left\{
\begin{array}{ll}
\dfrac{1}{8}\left( b-a\right) \left\| f^{\prime }\right\| _{\left[ a,b\right]
,\infty } & \text{if \hspace{0.05in}}f^{\prime }\in L_{\infty }\left[ a,b%
\right] ; \\
&  \\
\dfrac{1}{4}\cdot \dfrac{\left( b-a\right) ^{\frac{1}{q}}}{\left( q+1\right)
^{\frac{1}{q}}}\left\| f^{\prime }\right\| _{\left[ a,b\right] ,p} & \text{%
if \hspace{0.05in}}f^{\prime }\in L_{p}\left[ a,b\right] ,\;p>1,\;\frac{1}{p}%
+\frac{1}{q}=1; \\
&  \\
\dfrac{1}{4}\left\| f^{\prime }\right\| _{\left[ a,b\right] ,1} & \text{if
\hspace{0.05in}}f^{\prime }\in L_{1}\left[ a,b\right] .%
\end{array}
\right.
\end{multline}
Let $I_{n}:a=x_{0}<x_{1}<\dots <x_{n-1}<x_{n}=b$ be a division of the
interval $\left[ a,b\right] $ and $h_{i}:=x_{i+1}-x_{i}$ $\left( i=0,\dots
,n-1\right) $ and $\nu \left( I_{n}\right) :=\max \left\{ h_{i}|i=0,\dots
,n-1\right\} .$

Consider the composite quadrature rule
\begin{equation}
Q_{n}\left( I_{n},f\right) :=\frac{1}{2}\sum_{i=0}^{n-1}\left[ f\left( \frac{%
3x_{i}+x_{i+1}}{4}\right) +f\left( \frac{x_{i}+3x_{i+1}}{4}\right) \right]
h_{i}.  \label{3.2}
\end{equation}
The following result holds.

\begin{theorem}
\label{t3.1}Let $f:\left[ a,b\right] \rightarrow \mathbb{R}$ be an
absolutely continuous function on $\left[ a,b\right] .$ Then we have
\begin{equation}
\int_{a}^{b}f\left( t\right) dt=Q_{n}\left( I_{n},f\right) +R_{n}\left(
I_{n},f\right) ,  \label{3.3}
\end{equation}
where $Q_{n}\left( I_{n},f\right) $ is defined by formula (\ref{3.2}), and
the remainder satisfies the estimates
\begin{equation}
\left| R_{n}\left( I_{n},f\right) \right| \leq \left\{
\begin{array}{ll}
\dfrac{1}{8}\left\| f^{\prime }\right\| _{\left[ a,b\right] ,\infty
}\sum\limits_{i=0}^{n-1}h_{i}^{2} & \text{if \hspace{0.05in}}f^{\prime }\in
L_{\infty }\left[ a,b\right] ; \\
&  \\
\dfrac{1}{4\left( q+1\right) ^{\frac{1}{q}}}\left\| f^{\prime }\right\| _{%
\left[ a,b\right] ,p}\left( \sum\limits_{i=0}^{n-1}h_{i}^{q+1}\right) ^{%
\frac{1}{q}} & \text{if \hspace{0.05in}}f^{\prime }\in L_{p}\left[ a,b\right]
,\; \\
& p>1,\;\frac{1}{p}+\frac{1}{q}=1; \\
\dfrac{1}{4}\left\| f^{\prime }\right\| _{\left[ a,b\right] ,1}\nu \left(
I_{n}\right) . &
\end{array}
\right.  \label{3.4}
\end{equation}
\end{theorem}

\begin{proof}
Applying inequality (\ref{3.1}) on the intervals $\left[ x_{i},x_{i+1}\right]
,$ we may state that
\begin{multline}
\left| \int_{x_{i}}^{x_{i+1}}f\left( t\right) dt-\frac{1}{2}\left[ f\left(
\frac{3x_{i}+x_{i+1}}{4}\right) +f\left( \frac{x_{i}+3x_{i+1}}{4}\right) %
\right] h_{i}\right|  \label{3.5} \\
\leq \left\{
\begin{array}{l}
\dfrac{1}{8}h_{i}^{2}\left\| f^{\prime }\right\| _{\left[ x_{i},x_{i+1}%
\right] ,\infty } \\
\\
\dfrac{1}{4\left( q+1\right) ^{\frac{1}{q}}}h_{i}^{1+\frac{1}{q}}\left\|
f^{\prime }\right\| _{\left[ x_{i},x_{i+1}\right] ,p},\;\;\;p>1,\;\frac{1}{p}%
+\frac{1}{q}=1; \\
\\
\dfrac{1}{4}h_{i}\left\| f^{\prime }\right\| _{\left[ x_{i},x_{i+1}\right]
,1};%
\end{array}
\right.
\end{multline}
for each $i\in \left\{ 0,\dots ,n-1\right\} .$

Summing the inequality (\ref{3.5}) over $i$ from $0$ to $n-1$ and using the
generalised triangle inequality, we get
\begin{equation}
\left| R_{n}\left( I_{n},f\right) \right| \leq \left\{
\begin{array}{l}
\dfrac{1}{8}\sum\limits_{i=0}^{n-1}h_{i}^{2}\left\| f^{\prime }\right\| _{%
\left[ x_{i},x_{i+1}\right] ,\infty } \\
\\
\dfrac{1}{4\left( q+1\right) ^{\frac{1}{q}}}\sum\limits_{i=0}^{n-1}h_{i}^{1+%
\frac{1}{q}}\left\| f^{\prime }\right\| _{\left[ x_{i},x_{i+1}\right]
,p},\;\;p>1,\;\frac{1}{p}+\frac{1}{q}=1; \\
\\
\dfrac{1}{4}\sum\limits_{i=0}^{n-1}h_{i}\left\| f^{\prime }\right\| _{\left[
x_{i},x_{i+1}\right] ,1}.%
\end{array}
\right.  \label{3.6}
\end{equation}
Now, we observe that
\begin{equation*}
\sum\limits_{i=0}^{n-1}h_{i}^{2}\left\| f^{\prime }\right\| _{\left[
x_{i},x_{i+1}\right] ,\infty }\leq \left\| f^{\prime }\right\| _{\left[ a,b%
\right] ,\infty }\sum\limits_{i=0}^{n-1}h_{i}^{2}.
\end{equation*}
Using H\"{o}lder's discrete inequality, we may write that
\begin{align*}
\sum\limits_{i=0}^{n-1}h_{i}^{1+\frac{1}{q}}\left\| f^{\prime }\right\| _{%
\left[ x_{i},x_{i+1}\right] ,p}& \leq \left(
\sum\limits_{i=0}^{n-1}h_{i}^{\left( 1+\frac{1}{q}\right) q}\right) ^{\frac{1%
}{q}}\left( \sum\limits_{i=0}^{n-1}\left\| f^{\prime }\right\| _{\left[
x_{i},x_{i+1}\right] ,p}^{p}dt\right) ^{\frac{1}{p}} \\
& =\left( \sum\limits_{i=0}^{n-1}h_{i}^{q+1}\right) ^{\frac{1}{q}}\left(
\sum\limits_{i=0}^{n-1}\int_{x_{i}}^{x_{i+1}}\left| f^{\prime }\left(
t\right) \right| ^{p}dt\right) ^{\frac{1}{p}} \\
& =\left( \sum\limits_{i=0}^{n-1}h_{i}^{q+1}\right) ^{\frac{1}{q}}\left\|
f^{\prime }\right\| _{\left[ a,b\right] ,p}.
\end{align*}
Also, we note that
\begin{align*}
\sum\limits_{i=0}^{n-1}h_{i}\left\| f^{\prime }\right\| _{\left[
x_{i},x_{i+1}\right] ,1}& \leq \max\limits_{i=\overline{0,n-1}}\left\{
h_{i}\right\} \sum\limits_{i=0}^{n-1}\left\| f^{\prime }\right\| _{\left[
x_{i},x_{i+1}\right] ,1} \\
& =\nu \left( I_{n}\right) \left\| f^{\prime }\right\| _{\left[ a,b\right]
,1}.
\end{align*}
Consequently, by the use of (\ref{3.6}), we deduce the desired result (\ref%
{3.4}).
\end{proof}

For the particular case where the division $I_{n}$ is equidistant, i.e.,
\begin{equation*}
I_{n}:x_{i}=a+i\cdot \frac{b-a}{n},\;\;i=0,\dots ,n,
\end{equation*}
we may consider the quadrature rule:
\begin{multline}
Q_{n}\left( f\right) :=\frac{b-a}{2n}\sum\limits_{i=0}^{n-1}\left\{ f\left[
a+\left( \frac{4i+1}{4n}\right) \left( b-a\right) \right] \right.
\label{3.7} \\
+\left. f\left[ a+\left( \frac{4i+3}{4n}\right) \left( b-a\right) \right]
\right\} .
\end{multline}

The following corollary will be more useful in practice.

\begin{corollary}
\label{c3.2}With the assumption of Theorem \ref{t3.1}, we have
\begin{equation}
\int_{a}^{b}f\left( t\right) dt=Q_{n}\left( f\right) +R_{n}\left( f\right) ,
\label{3.8}
\end{equation}
where $Q_{n}\left( f\right) $ is defined by (\ref{3.7}) and the remainder $%
R_{n}\left( f\right) $ satisfies the estimate:
\begin{equation}
\left| R_{n}\left( I_{n},f\right) \right| \leq \left\{
\begin{array}{l}
\dfrac{1}{8}\left\| f^{\prime }\right\| _{\left[ a,b\right] ,\infty }\dfrac{%
\left( b-a\right) ^{2}}{n} \\
\\
\dfrac{1}{4\left( q+1\right) ^{\frac{1}{q}}}\left\| f^{\prime }\right\| _{%
\left[ a,b\right] ,p}\cdot \dfrac{\left( b-a\right) ^{1+\frac{1}{q}}}{n} \\
\\
\dfrac{1}{4}\left\| f^{\prime }\right\| _{\left[ a,b\right] ,1}\cdot \dfrac{%
\left( b-a\right) }{n}.%
\end{array}
\right.  \label{3.9}
\end{equation}
\end{corollary}

\section{Applications for P.D.F.'s}

Summarising some of the results in Section 2, we may state that for $f:\left[
a,b\right] \rightarrow \mathbb{R}$ an absolutely continuous function, we
have the inequality
\begin{multline}
\left| \frac{1}{2}\left[ g\left( x\right) +g\left( a+b-x\right) \right] -%
\frac{1}{b-a}\int_{a}^{b}f\left( t\right) dt\right|  \label{4.1} \\
\leq \left\{
\begin{array}{l}
\left[ \dfrac{1}{8}+2\left( \dfrac{x-\frac{3a+b}{4}}{b-a}\right) ^{2}\right]
\left( b-a\right) \left\| g^{\prime }\right\| _{\left[ a,b\right] ,\infty
}\hfill \text{if }g^{\prime }\in L_{\infty }\left[ a,b\right] \\
\\
\dfrac{2^{\frac{1}{q}}}{\left( q+1\right) ^{\frac{1}{q}}}\left[ \left(
\dfrac{x-a}{b-a}\right) ^{q+1}+\left( \dfrac{\frac{a+b}{2}-x}{b-a}\right)
^{q+1}\right] ^{\frac{1}{q}}\left( b-a\right) ^{\frac{1}{q}}\left\|
g^{\prime }\right\| _{\left[ a,b\right] ,p}, \\
\hfill \text{if \hspace{0.05in}}p>1,\;\frac{1}{p}+\frac{1}{q}=1,\text{ and }%
g^{\prime }\in L_{p}\left[ a,b\right] ; \\
\left[ \dfrac{1}{4}+\left| \dfrac{x-\frac{3a+b}{4}}{b-a}\right| \right]
\left\| g^{\prime }\right\| _{\left[ a,b\right] ,1},%
\end{array}
\right.
\end{multline}
for all $x\in \left[ a,\frac{a+b}{2}\right] .$

Now, let $X$ be a random variable taking values in the finite interval $%
\left[ a,b\right] ,$ with the probability density function $f:\left[ a,b%
\right] \rightarrow \lbrack 0,\infty )$ and with the cumulative distribution
function $F\left( x\right) =\Pr \left( X\leq x\right) =\int_{a}^{x}f\left(
t\right) dt.$

The following result holds.

\begin{theorem}
\label{t4.1}With the above assumptions, we have the inequality
\begin{multline}
\left| \frac{1}{2}\left[ F\left( x\right) +F\left( a+b-x\right) \right] -%
\frac{b-E\left( X\right) }{b-a}\right|  \label{4.2} \\
\leq \left\{
\begin{array}{l}
\left[ \dfrac{1}{8}+2\left( \dfrac{x-\frac{3a+b}{4}}{b-a}\right) ^{2}\right]
\left( b-a\right) \left\| f\right\| _{\left[ a,b\right] ,\infty }\hfill
\text{if }f\in L_{\infty }\left[ a,b\right] \\
\\
\dfrac{2^{\frac{1}{q}}}{\left( q+1\right) ^{\frac{1}{q}}}\left[ \left(
\dfrac{x-a}{b-a}\right) ^{q+1}+\left( \dfrac{\frac{a+b}{2}-x}{b-a}\right)
^{q+1}\right] ^{\frac{1}{q}}\left( b-a\right) ^{\frac{1}{q}}\left\|
f\right\| _{\left[ a,b\right] ,p}, \\
\hfill \text{if \hspace{0.05in}}p>1,\;\frac{1}{p}+\frac{1}{q}=1,\text{ and }%
f\in L_{p}\left[ a,b\right] ; \\
\left[ \dfrac{1}{4}+\left| \dfrac{x-\frac{3a+b}{4}}{b-a}\right| \right] ,%
\end{array}
\right.
\end{multline}
for any $x\in \left[ a,\frac{a+b}{2}\right] ,$ where $E\left( X\right) $ is
the expectation of $X.$
\end{theorem}

\begin{proof}
Follows by (\ref{4.1}) on choosing $g=F$ and taking into account that
\begin{equation*}
E\left( X\right) =\int_{a}^{b}tdF\left( t\right) =b-\int_{a}^{b}F\left(
t\right) dt.
\end{equation*}
\end{proof}

In particular, we have:

\begin{corollary}
\label{cnew}With the above assumptions, we have
\begin{multline}
\left| \frac{1}{2}\left[ F\left( \frac{3a+b}{4}\right) +F\left( \frac{a+3b}{4%
}\right) \right] -\frac{b-E\left( X\right) }{b-a}\right|  \label{4.3} \\
\leq \left\{
\begin{array}{l}
\dfrac{1}{8}\left( b-a\right) \left\| f\right\| _{\left[ a,b\right] ,\infty
}\hfill \text{if }f\in L_{\infty }\left[ a,b\right] \\
\\
\dfrac{1}{4}\cdot \dfrac{\left( b-a\right) ^{\frac{1}{q}}}{\left( q+1\right)
^{\frac{1}{q}}}\left\| f\right\| _{\left[ a,b\right] ,p},\;\hfill \text{if
\hspace{0.05in}}p>1,\;\frac{1}{p}+\frac{1}{q}=1,\text{ and }f\in L_{p}\left[
a,b\right] ; \\
\\
\dfrac{1}{4}.%
\end{array}
\right.
\end{multline}
\end{corollary}

\end{document}